\documentclass[preprint]{elsarticle}

\usepackage[bitstream-charter]{mathdesign}

\usepackage{latexsym}
\usepackage{amsmath}
\usepackage{amsthm}
\usepackage{graphicx}
\usepackage{enumerate}

\newtheorem{theorem}{Theorem}[section]
\newtheorem{lemma}[theorem]{Lemma}
\newtheorem{corollary}[theorem]{Corollary}
\newtheorem{proposition}[theorem]{Proposition}

\newtheoremstyle{publishedexample}
  {\topsep}{\topsep}{\normalfont}{}{\itshape}{.}{.5em}{}
\theoremstyle{publishedexample}
\newtheorem{example}[theorem]{Example}
\theoremstyle{plain}

\newproof{pf}{Proof}

\numberwithin{equation}{section}

\makeatletter
\renewcommand{\@biblabel}[1]{#1.}
\makeatother

\begin{document}

\pagestyle{plain}

\noindent
\begin{minipage}[t]{0.58\textwidth}
\small
Ann. Comb. 19 (2015) 125--130\\
DOI 10.1007/s00026-015-0259-8\\
Published online January 11, 2015\\
\copyright\ Springer Basel 2015
\end{minipage}%
\begin{minipage}[t]{0.42\textwidth}
\raggedleft\textbf{Annals of Combinatorics}
\end{minipage}

\bigskip\bigskip
\noindent{\LARGE\bfseries A Characterization of Tangle Matroids\par}
\bigskip
\noindent Dennis Hall\par
\smallskip
\noindent Department of Mathematics, Louisiana State University, Baton Rouge, Louisiana\par
\noindent dhall15@math.lsu.edu\par
\medskip
\noindent Received April 25, 2013\par
\medskip
\noindent\hfill\textit{Mathematics Subject Classification:} 05B35\par

\begin{center}
Dedicated to James Oxley, for his years of patience and guidance
\end{center}

\medskip
\noindent\textbf{Abstract.} We prove that a matroid with at least two elements is a tangle matroid if and only if it cannot be covered by three hyperplanes. Some consequences of this theorem are also noted. In particular, no binary matroid of rank at least two is a tangle matroid.


\bigskip

\section{Introduction}

Tangles are a way of identifying highly connected regions in matroids of low connectivity. As these highly-connected regions are the obstructions to large branchwidth (see \cite{GeelenObstruction}), tangles play a fundamental role in Geelen, Gerards, and Whittle's structure theorem for $GF(q)$-representable matroids in \cite{structureTheorem}. The reader is also referred to \cite{structureTheorem} for a history and enlightening discussion of tangles.

The matroid terminology used here will follow Oxley \cite{oxley2011matroid}. In particular, the connectivity function of a matroid $M$ is defined, for $X\subseteq E(M)$, to be \[\lambda_M(X)=r(X)+r(E(M)-X)-r(M).\] Our definitions of tangles and of tangle matroids follow \cite{Stefan}. Let $M$ be a matroid, and $\mathcal{T}$ a collection of subsets of $E(M)$. Then $\mathcal{T}$ is a {\it tangle of order $\theta$ of $M$} if
\begin{enumerate}
\item [(T1)]for all $X\in\mathcal{T}$, $\lambda_M(X)<\theta$;
\item [(T2)]for all $X\subseteq E(M)$ with $\lambda_M(X)<\theta$, either $X\in\mathcal{T}$ or $E(M)-X\in\mathcal{T}$;
\item [(T3)]if $X,Y,Z\in\mathcal{T}$, then $X\cup Y\cup Z\neq E(M)$;
\item [(T4)]for each $e\in E(M)$, $E(M)-\{e\}\notin\mathcal{T}$.
\end{enumerate}

A useful tool for studying tangles is the {\it tangle matroid}. Let $M$ be a matroid and $\mathcal{T}$ a tangle of $M$ of order $\theta$. Let $\rho:2^{E(M)}\rightarrow\mathbb{N}$ be defined by
	\begin{align*}
		\rho(X) := \left\{\begin{array}{ll}
			\min \{ \lambda_M(Y) : X \subseteq Y \in \mathcal{T} \}, & \textrm{ if there exists a $Y$ with } X \subseteq Y \in \mathcal{T},\\   
			\theta, & \textrm{ otherwise}.
		\end{array}\right.
	\end{align*}

In \cite{GeelenObstruction}, it is shown that $\rho$ is the rank function of a matroid on $E(M)$. We call this matroid $M(\mathcal{T})$. If $\mathcal{T}$ is a tangle of a matroid $M$, we call $M(\mathcal{T})$ a {\it tangle matroid of} $M$.

In this paper, we characterize those matroids that are tangle matroids. For example, $PG(2,3)$ is a tangle matroid while $PG(2,2)$ is not. The reason lies in the arrangement of the hyperplanes in these matroids. In Section \ref{secTM}, we prove the following result.

\begin{theorem}\label{main}
A matroid $M$ other than $U_{1,1}$ is a tangle matroid if and only if $M$ has no three hyperplanes whose union is $E(M)$.
\end{theorem}

In addition, we show some consequences of this theorem. In particular, a tangle matroid of $M$ is a quotient of $M$, and all binary tangle matroids have rank less than two.

\section{High-Order Tangles}

An immediate consequence of the following lemma \cite{geelen} is that a matroid cannot have a tangle whose order exceeds its rank.

\begin{lemma}\label{geelenLemma}
Let $\mathcal{T}$ be a tangle of order $\theta$ at least two in a matroid $M$. Then each subset of $E(M)$ with rank less than $\theta$ is in $\mathcal{T}$.
\end{lemma}

In view of this, it is natural to ask when a matroid has a tangle whose order equals the rank of the matroid.

\begin{proposition}\label{hyperplanes}
A matroid $M$ other than $U_{1,1}$ has an order-$r(M)$ tangle if and only if $M$ has no three hyperplanes whose union is $E(M)$.
\end{proposition}
\begin{proof}
If a matroid has three such hyperplanes, then, from Lemma $\ref{geelenLemma}$, each of these hyperplanes is a member of the tangle, a contradiction to (T3). On the other hand, if no three hyperplanes cover the entire ground set of $M$, we show that
\[
\mathcal{T}=\{X:r(X)<r(M)\}
\]
is an order-$r(M)$ tangle of $M$.

For (T1), we choose $X\in\mathcal{T}$. If $r(E(M)-X)<r(M)$, then $X$ and $E(M)-X$ are two non-spanning sets that cover $E(M)$. As such, there are two hyperplanes $H_1$ and $H_2$ containing $X$ and $E(M)-X$, respectively, such that $H_1\cup H_2=E(M)$. Thus $r(E(M)-X)=r(M)$, and therefore, $\lambda_M(X)=r(X)<r(M)$, as desired.

For (T2), choose a partition $(X,Y)$ of $E(M)$ such that $\lambda_M(X)<r(M)$. Expanding out $\lambda_M(X)$, this becomes $r(X)+r(Y)<2r(M)$. It follows, then, that either $r(X)<r(M)$ or $r(Y)<r(M)$, and thus either $X\in\mathcal{T}$ or $Y\in\mathcal{T}$.

The remaining conditions, (T3) and (T4), follow immediately from the hypothesis. We conclude that $\mathcal{T}$ is indeed a tangle of order $r(M)$.
\end{proof}

We now know which matroids have a tangle with order equal to the rank of the matroid. A natural next challenge is to find all such tangles. It turns out that, if a matroid $M$ has an order-$r(M)$ tangle, then that tangle is exactly the one used in the proof of Proposition \ref{hyperplanes}. The following proposition shows this.

\begin{proposition}\label{hyperplanes2}
If a matroid $M$ has an order-$r(M)$ tangle $\mathcal{T}$, then $\mathcal{T}=\{X:r(X)<r(M)\}$.
\end{proposition}
\begin{proof}
Set $\mathcal{T}^\prime=\{X:r(X)<r(M)\}$ and let $\mathcal{T}$ be an order-$r(M)$ tangle of $M$. Assume first that $r(M)=1$. If $x$ is a loop of $M$, then, by combining (T2) and (T3), we have that $\{x\}\in\mathcal{T}$. From this and \cite[Lemma 2.27]{Stefan}, it follows that, if $X\subseteq E(M)$ where $r(X)=0$, then $X\in\mathcal{T}$. Thus $\mathcal{T}^\prime\subseteq\mathcal{T}$ and we have only to show that $\mathcal{T}$ contains no rank-$1$ elements. Assuming the contrary, let $A\subseteq E(M)$ be such that $r(A)=1$ and $A\in\mathcal{T}$. Since, in this case, $\lambda_M(A)=0$, it follows that $r(E(M)-A)=0$. From the above argument, $E(M)-A\in\mathcal{T}$, contradicting (T3) since both $A$ and its complement would be in $\mathcal{T}$. Therefore, $\mathcal{T}^\prime=\mathcal{T}$.

If $r(M)\geq 2$, then we are able to use Lemma \ref{geelenLemma} to show that $\mathcal{T}^\prime\subseteq\mathcal{T}$. To prove the equality, suppose $Y\in\mathcal{T}$ such that $r(Y)=r(M)$. Then, by (T1), \[r(M)>\lambda_M(Y)=r(Y)+r(E(M)-Y)-r(M)=r(E(M)-Y).\] Thus, since $r(E(M)-Y)<r(M)$, both $Y$ and its complement are members of $\mathcal{T}$, contradicting (T3). Therefore, $\mathcal{T}^\prime=\mathcal{T}$.
\end{proof}

\section{Tangle Matroids}\label{secTM}

The first question we answer concerning tangle matroids is that of when a matroid has itself as a tangle matroid.

\begin{proposition}\label{tangleMatroids1}A matroid $M$ other than $U_{1,1}$ has itself as a tangle matroid if and only if $M$ has no three hyperplanes whose union is $E(M)$.\end{proposition}
\begin{proof}
Let $\mathcal{T}$ be a tangle of $M$ for which $M(\mathcal{T})=M$. Then the order of $\mathcal{T}$ is necessarily $r(M)$ and thus, from Proposition \ref{hyperplanes}, $M$ has no three hyperplanes whose union is $E(M)$. We thus assume that $M$ has no three hyperplanes whose union is $E(M)$. Then, by Proposition \ref{hyperplanes}, $M$ has an order-$r(M)$ tangle $\mathcal{T}$. Further, Proposition \ref{hyperplanes2} tells us that $\mathcal{T}=\{X:r(X)<r(M)\}$. We now have only to show that $M(\mathcal{T})=M$.

Choose $X\subseteq E(M)$. We wish to show that $\rho(X)=r(X)$. If $r(X)=r(M)$, then neither $X$, nor any set containing $X$, is a member of $\mathcal{T}$. Thus $\rho(X)=\theta=r(M)$. On the other hand, if $r(X)<r(M)$, then $X\in\mathcal{T}$ and $\lambda_M(X)=r(X)+r(E(M)-X)-r(M)$. Since $\mathcal{T}$ is, by Proposition \ref{hyperplanes2}, a tangle, we use (T3) to observe that $E-X\not\in\mathcal{T}$. Thus $r(E-X)=r(M)$ and therefore $\lambda_M(X)=r(X)$. If $X\subseteq Y\in\mathcal{T}$, then, as above we see that $\lambda_M(Y)=r(Y)\geq r(X)=\lambda_M(X)$. Thus,
\[
\rho(X)=\min \{ \lambda_M(Y) : X \subseteq Y \in \mathcal{T} \} = \lambda_M(X)=r(X),
\]
and therefore, $M(\mathcal{T})=M$.
\end{proof}

We are now ready to prove the main theorem. We show that a matroid $M$ other than $U_{1,1}$ is a tangle matroid if and only if it has no three hyperplanes that cover $E(M)$. We omit $U_{1,1}$ in this statement as it is the only matroid without three covering hyperplanes that is not a tangle matroid.

\begin{proof}[Proof of Theorem 1.1]
If $M$ has no three hyperplanes whose union is $E(M)$, then, by Proposition \ref{tangleMatroids1}, $M$ has a tangle $\mathcal{T}$ such that $M=M(\mathcal{T})$. We thus assume $M$ has three hyperplanes $H_1$, $H_2$, and $H_3$ such that $H_1\cup H_2\cup H_3=E(M)$. Assume that $\mathcal{T}$ is a tangle for which $M(\mathcal{T})=M$. Observe that the order of $\mathcal{T}$ is $r(M)$. Thus, from the definition of $M(\mathcal{T})$, there exist sets $Y_1$, $Y_2$, and $Y_3$ such that each $H_i\subseteq Y_i\in\mathcal{T}$. It follows, then, that $Y_1\cup Y_2\cup Y_3=E(M)$, contradicting (T3).
\end{proof}

We now list several immediate consequences of Theorem \ref{main}. The first gives several ways to restate the main theorem.

\begin{corollary}
Let $M$ be a matroid with ground set $E$. Then the following are equivalent.
\begin{enumerate}
\item $M$ is either a tangle matroid or $U_{1,1}$.
\item $M$ has no three hyperplanes whose union is $E$.
\item In every partition $(X,Y,Z)$ of $E$, at least one of $X,Y,$ and $Z$ is a spanning set.
\item $M$ has no three cocircuits whose intersection is empty.
\end{enumerate}
\end{corollary}

A matroid is {\it round} if it has no two disjoint cocircuits. The following is another immediate consequence of Theorem \ref{main}.

\begin{corollary}
Every tangle matroid is round.
\end{corollary}

The next consequence follows immediately from Theorem \ref{main} and Geelen and Kabell's result \cite{EPMatroids} that finds an Erd{\H o}s-P\'osa property for matroid circuits. The bicircular matroid of a graph $G$ is denoted by $B(G)$.

\begin{corollary}
There exists an integer-valued function $\gamma(n)$ such that, for any $n\in\mathbb{N}$, if $M$ is a tangle matroid with $r(M)\geq\gamma(n)$, then $M$ has a minor isomorphic to $U_{n,2n}$, $M(K_n)$, or $B(K_n)$.
\end{corollary}

Observe that the class of tangle matroids is not closed under minors. We can take, for example, any single-element deletion of $U_{2,4}$. The problem with minors arises from deleting elements that are not contained in a parallel class. If we do not allow this type of deletion, that is, if we look at parallel minors, then this problem is eliminated. The following corollary states this result.

\begin{corollary}\label{parmincor}
If $M$ is a tangle matroid, then every parallel minor of $M$ is either a tangle matroid or $U_{1,1}$.
\end{corollary}

Let $\mathcal{C}_\mathcal{T}$ be the class consisting of $U_{1,1}$ together with all tangle matroids. Corollary \ref{parmincor} establishes that $\mathcal{C}_\mathcal{T}$ is closed under taking parallel minors. One might hope, then, to obtain an excluded parallel minor-based characterization of $\mathcal{C}_\mathcal{T}$. Unfortunately, such a list of excluded minors would be infinite. For example, the class $\{U_{k,3k-3}:k\geq 2\}$ are all excluded parallel minors of $\mathcal{C}_\mathcal{T}$.

If we restrict $\mathcal{C}_\mathcal{T}$ to only its $GF(q)$-representable members, a similar problem occurs. Let $M_q$ be the class of $GF(q)$-representable matroids. If $q\geq 3$, then the class $M_q\cap \mathcal{C}_\mathcal{T}$ has an infinite number of excluded parallel minors, as we show in the following example.

\begin{example}
Choose two integers $r,q\geq 3$. Let $M_1$ and $M_2$ be two matroids that are isomorphic to $PG(r-1,q)$ such that $E(M_1)\cap E(M_2)=T$ where $M_1|T=M_2|T\cong PG(r-2,q)$. Recall that $P_T(M_1,M_2)$ denotes the generalized parallel connection of $M_1$ and $M_2$ over $T$ (see \cite{oxley2011matroid}). Then $M=P_T(M_1,M_2)\backslash T$ is not a tangle matroid since $E(M_1)-T$ and $E(M_2)-T$ are two hyperplanes whose union covers $E(M)$. However, for $a\in E(M)$, $M/a$ is isomorphic to $M_1$, where each member of $T$ is replaced with $q-1$ elements in parallel. Thus $M/a$ is a $GF(q)$-representable tangle matroid and so $M$ is an excluded parallel minor of $M_q\cap \mathcal{C}_\mathcal{T}$.
\end{example}

We conclude this section by characterizing the binary tangle matroids.

\begin{theorem}\label{thm:binary}
A matroid $M$ other than $U_{1,1}$ is a binary tangle matroid if and only if $r(M)<2$.
\end{theorem}
\begin{proof}
Let $M$ be a binary tangle matroid. It is well known (see, for example, \cite{oxley2011matroid}), that a matroid is binary if and only if the matroid has no parallel minor isomorphic to $U_{2,k+2}$ for any $k\geq 2$.  Suppose $r(M)\geq 2$. Let $H$ be an independent set of $M$ such that $|H|=r(M)-2$. Then $si(M/H)$ is a rank-$2$ parallel minor of $M$ and, as both $\mathcal{C_\mathcal{T}}$ and $M_2$ are closed under parallel minors, $si(M/H)$ is a member of $\mathcal{C_\mathcal{T}}\cap M_2$. However, since neither $U_{2,2}$ nor $U_{2,3}$ is a tangle matroid, $si(M/H)$ is isomorphic to one of the excluded parallel minors for binary matroids, a contradiction. Therefore, $r(M)<2$. The converse follows immediately from Theorem \ref{main}.
\end{proof}

\section{Quotients}

It is easy to check that, if a matroid $M$ has a tangle matroid $M^\prime$, then $r(M^\prime)\leq r(M)$. Because of this, a tangle matroid $M^\prime$ can be thought of as a minimal matroid with $M^\prime$ as a tangle matroid. The following theorem is used, by Corollary \ref{TangQuot}, to show this formally, where the ordering is that induced by matroid quotients.

\begin{theorem}\label{closure}
Let $M$ be a matroid and $X\subseteq E(M)$. If $\mathcal{T}$ is a tangle of $M$, then $cl_M(X)\subseteq cl_\mathcal{T}(X)$.
\end{theorem}
\begin{proof}
Assume $\theta\geq 1$ and $|E(M)| > 3$. If $\rho(X)=\theta$ then the result is trivial. Pick $X$ with $\rho(X)=t<\theta$, and $a\in cl_M(X)-X$. There exists a set $Y\in\mathcal{T}$ with $X\subseteq Y$ and $\lambda_M(Y)=t$. Since $a\in cl_M(Y)$, we have $\lambda_M(Y\cup a)\leq\lambda_M(Y)=t$. Since $Y,\{a\}\in\mathcal{T}$, also $Y\cup a\in\mathcal{T}$. Since $X\subseteq Y\cup a$, we have $\lambda_M(Z)\geq t$ for all $Z\in\mathcal{T}$ with $Y\cup a\subseteq Z$, and therefore, $\rho(Y\cup a)=t$. Now $\rho(X)=\rho(Y\cup a)$ and therefore, $Y\cup a\subseteq cl_\mathcal{T}(X)$.
\end{proof}

The following corollary follows immediately from Theorem \ref{closure} by using results from \cite[\S7.3]{oxley2011matroid}.

\begin{corollary}\label{TangQuot}
Let $M$ be a matroid with tangle $\mathcal{T}$. Then the following statements hold.
\begin{enumerate}[(i)]
\item $M(\mathcal{T})$ is a quotient of $M$.
\item Every flat of $M(\mathcal{T})$ is a flat of $M$.
\item If $X\subseteq Y\subseteq E(M)$, then $r_M(Y)-r_M(X)\geq r_\mathcal{T}(Y)-r_\mathcal{T}(X)$.
\item Every circuit of $M$ is a union of circuits of $M(\mathcal{T})$.
\item The identity map $i:M\longrightarrow M(\mathcal{T})$ is a strong map.
\end{enumerate}
\end{corollary}

The result that $M(\mathcal{T})$ is a quotient of $M$ has the following consequence.

\begin{proposition}
Let $M$ be a matroid such that both $M$ and $M^*$ are tangle matroids. Then $M=M^*$.
\end{proposition}
\begin{proof}
If $M$ is a tangle matroid, then there is a tangle $\mathcal{T}$ of $M$ such that $M=M(\mathcal{T})$. As $M$ and $M^*$ share the same connectivity function $\lambda_M$, it follows that $\mathcal{T}$ is a tangle of $M^*$ and thus $M(\mathcal{T})=M$ is a quotient of $M^*$. Similarly, $M^*$ is a quotient of $M$. Therefore, $M=M^*$.
\end{proof}

\par\medskip\noindent\textbf{Acknowledgments.}
The author thanks James Oxley and an anonymous referee for their valuable advice in the preparation of this paper.

\end{document}